\documentclass[11pt]{article}
\usepackage[margin=1in]{geometry}
\usepackage{graphicx} 
\usepackage{hyperref}   
\hypersetup{
	colorlinks=true,
	linkcolor=blue,
	citecolor=blue,
	urlcolor=blue,
}
\usepackage[square,numbers]{natbib}
 
\usepackage{graphicx}
\usepackage{subcaption}
\usepackage{amsmath,mathtools,amssymb,amsfonts}
\usepackage{xcolor}
\usepackage{algorithm}
\usepackage{algpseudocode}
\usepackage{nicematrix}

\usepackage{mathtools,amssymb,lipsum}

\usepackage{cuted}
\newcommand{\norm}[1]{\left\lVert#1\right\rVert}
\newcommand{\xddots}{%
  \raise 4pt \hbox {.}
  \mkern 6mu
  \raise 1pt \hbox {.}
  \mkern 6mu
  \raise -2pt \hbox {.}
}

\newtheorem{definition}{Definition}

\title{Functional Subspace Projection for Detection of Coordinated Stealthy Attacks in Power Systems 
}

\title{Functional Subspace Projection for Detection of Coordinated
Stealthy Attacks in Power Systems}
\author{Koto Omiloli, Satish Vedula and Olugbenga Moses Anubi}

\date{}

\begin{document}
\maketitle
\begin{center}
Department of Electrical and Computer Engineering\\
Center for Advanced Power Systems, Florida State University\\
E-mail: \{kao23a, svedula, oanubi\}@fsu.edu
\end{center} 

\begin{center}
\small
\textit{A revised version of this manuscript has been accepted for publication in
ASME Letters in Dynamic Systems and Control (ALDSC).}
\end{center}

\vspace{0.5em}
\begin{abstract}
{\it \textbf{Abstract:} Modern power grids are increasingly vulnerable to coordinated cyber-attacks, particularly false data injection attacks  (FDIAs) that can evade conventional residual-based detectors. While most existing detection methods rely on instantaneous measurements, coordinated dynamic attacks can remain stealthy at each time step while introducing structured temporal deviations. This paper develops a joint framework for modeling and detecting such attacks in multi-area power systems. A time-aggregated attack model is first formulated to capture temporal evolution and inter-area coordination. For detection, a kernel-embedded functional subspace detection (KEFSD) method is proposed, which models residual trajectories in a reproducing kernel Hilbert space (RKHS) and employs RKHS-constrained functional principal component analysis (PCA) to identify anomalous temporal patterns. Simulation results on a modified IEEE 14-bus system demonstrates the proposed method achieves improved detection performance compared to the conventional residual based 2-norm detector. 
}
\end{abstract}


\section{Introduction}
The increasing integration of communication networks, distributed generation, and advanced control architectures has transformed modern power grids into tightly coupled cyber–physical systems. While these capabilities enhance operational flexibility and efficiency, they also expose the grid to sophisticated cyber-attacks targeting measurement and control channels. Among these threats, false data injection attacks (FDIAs) are particularly concerning because carefully designed perturbations can evade classical bad data detection mechanisms, such as residual-based and weighted least squares (WLS) detectors embedded within state estimation algorithms.

Beyond centralized estimation, vulnerabilities extend to distributed architectures as well. For example,~\cite{xu2024globally} demonstrated that coordinated attacks can bias distributed state estimators without triggering alarms, showing that stealthiness can persist under distributed consistency requirements. \cite{liu2020network} proposed coordinated attacks that jointly manipulate network parameters and system states, reducing the number of measurements that must be compromised. In integrated infrastructures, ~\cite{zadsar2022prevention} formulated bilevel optimization–based attacks on coupled power–gas systems, capable of bypassing detection while disrupting gas fired generation. Collectively, these studies highlight the increasing sophistication of multi-node and multi-infrastructure attacks.

To mitigate such threats, data-driven detection and resilient estimation methods have been explored. \cite{ashok2016online} proposed an online framework using load forecasts, generation schedules, and phase-measurement-units (PMU) to detect stealthy anomalies beyond traditional residual checks, emphasizing adaptive real-time detection. Its effectiveness, however depends on forecast quality and PMU coverage. Similarly, \cite{chen2022data} developed a deep autoencoding Gaussian mixture model for unsupervised detection, which is effective statistically but may fail if attackers mimic the learned data distribution. Supervised, semi-supervised, and online learning approaches \cite{ozay2015machine, wang2019detection} have also been used to identify stealthy attacks, yet they do not explicitly account for physical system constraints, leaving some attacks undetected. For resilient estimation, \cite{anubi2019enhanced} proposed a method leveraging Gaussian process priors within convex optimization frameworks to recover system states under corrupted measurements, but offer limited adaptability to coordinated, multi-area attacks.

In parallel, efforts to improve WLS-based detection have been considered, including robust estimators such as least winsorized squares~\cite{basetti2015robust} and ensemble anomaly detection integrated with WLS~\cite{gholami2024detection}. While these methods enhance resilience against gross errors and outliers, they remain fundamentally residual-based and cannot eliminate the geometric vulnerability that enables stealth attacks within the measurement subspace. To overcome this limitation, data-driven kernel and latent-space methods have gained attention. For example, \cite{zhang2024detection} applied kernel principal component analysis (KPCA) to detect AC-model FDIAs using high-dimensional statistical metrics; \cite{ma2020network} employed radial basis function (RBF) KPCA for nonlinear feature extraction followed by decision tree classification, and \cite{bhattacharjee2022deep} used self-supervised latent space clustering with autoencoders to model nonlinear measurement manifolds. These approaches improve statistical separability, though their effectiveness depends on the quality of the learned data manifold and does not inherently guarantee detection of attacks that respect the physical measurement relationships of the system.

Building on these insights, this paper develops a joint framework for modeling and detecting coordinated dynamic stealthy attacks in power systems. We formulate a time-aggregated multi-area attack model that captures inter-area correlations via graph-structured coupling. For detection, we propose a kernel-embedded functional subspace detection (KEFSD) method that models residual trajectories as smooth functions in a reproducing kernel Hilbert space (RKHS), enabling the identification of structured temporal deviations that remain undetected by conventional residual-based methods. To our knowledge, this is the first framework that jointly addresses dynamic physically consistent attack synthesis and RKHS-constrained functional detection.

The main contributions of this work are summarized as follows:
\begin{itemize}
\item We develop a dynamic multi-area FDIA model capturing temporal evolution across the power system measurements.
\item We derive a tractable stealthy attack design problem under linearized residual and sparsity constraints.
\item We propose an RKHS-based functional subspace detector using online projection of residual trajectories.
\item We validate the approach on a partitioned IEEE 14-bus system, showing improved detection over the conventional residual based 2-norm detector.
\end{itemize}

\section{Notation and Preliminaries}
We adopt the following notational conventions: \( \mathbb{R}^n \) denotes the space of real-valued vectors of dimension \( n \). Scalars are denoted by lowercase letters (e.g., $x \in \mathbb{R}$), and vectors by bold lowercase letters (e.g., $\mathbf{x} \in \mathbb{R}^n$), with entries denoted by $x_i$. Matrices are denoted by uppercase letters (e.g., \( K \in \mathbb{R}^{n \times m} \) with entries \( K_{ij}\) and column vectors denoted by $K[i]$.  The trace of a square matrix A is denoted by Tr(A). The zero vector is written as \( \mathbf{0}\). The identity matrix of size $n$ as \( I_n \).
A weighted undirected graph is denoted by $\mathcal{G} = (\mathcal{N}, \mathcal{E})$, 
where $\mathcal{N}$ is the set of nodes and $\mathcal{E}$ is the set of edges. 
The weighted adjacency matrix is denoted by $[w_{kl}] \in \mathbb{R}^{|\mathcal{N}| \times |\mathcal{N}|}$, 
where $w_{kl}$ represents the weight of the edge between nodes $k$ and $l$, and $w_{kl} = 0$ if $(k,l) \notin \mathcal{E}$. 

\begin{definition}
  A Hilbert space $\mathcal{H}$ of functions $f:\mathcal{X}\to\mathbb{R}$ is called a RKHS if, for every $x \in \mathcal{X}$, the point evaluation map $f \mapsto f(x)$ is a bounded linear functional on $\mathcal{H}$. Then there exists a function $\kappa:\mathcal{X}\times\mathcal{X}\to\mathbb{R}$ such that
$
f(x) = \langle f, \kappa(x,\cdot) \rangle_{\mathcal{H}}, \quad \forall f \in \mathcal{H},
$
where $\kappa$ is called the reproducing kernel.
\end{definition}

\vspace{-1em}
\section{Dynamic Power System Modeling}\label{dynamic_power_model}
We begin by representing the power network as a graph $\mathcal{G} = (\mathcal{N}, \mathcal{E})$, 
where $\mathcal{N} = \{1,\dots,N\}$ denotes the set of buses and 
$\mathcal{E}$ the set of transmission lines. Let 
$\mathcal{N}_G \subseteq \mathcal{N}$ and $\mathcal{N}_L \subseteq \mathcal{N}$ 
denote the sets of generator and load buses, respectively. Each synchronous generator $i \in \mathcal{N}_G$ is modeled using the classical
swing equations:

\vspace{-2em}
\begin{align}
\dot{\delta}_i &= \omega_i, \label{eq:swing1} \\
m_i \dot{\omega}_i &= p_{m,i} - p_{e,i} - d_i \omega_i, \label{eq:swing2}
\end{align}
\vspace{-2em}

where $\delta_i$ is the rotor angle, $\omega_i$ is the rotor frequency,
$m_i$ is the inertia constant, $d_i$ is the damping coefficient,
$p_{m,i}$ is the mechanical input power, and $p_{e,i}$ is the electrical output
power. The network algebraic constraints is given by:

\vspace{-2em}
\begin{align}
0 &= p_{e,i} - p_{l,i}
- \sum_{j \in \mathcal{N}} v_i v_j \Big(
G_{ij} \cos(\theta_i - \theta_j) \notag \\
&\qquad\qquad\quad
+ B_{ij} \sin(\theta_i - \theta_j)
\Big),
\label{eq:active_balance}
\\
0 &= q_{e,i} - q_{l,i}
- \sum_{j \in \mathcal{N}} v_i v_j \Big(
G_{ij} \sin(\theta_i - \theta_j) \notag \\
&\qquad\qquad\quad
- B_{ij} \cos(\theta_i - \theta_j)
\Big),
\label{eq:reactive_balance}
\end{align}
\vspace{-2em}

where $v_i, v_j$ are bus voltage magnitudes, $\theta_i, \theta_j$ are bus voltage angles, 
$G_{ij}$ and $B_{ij}$ are the real and imaginary parts of the bus admittance matrix $Y_\text{bus}$, 
respectively, $p_{e,i}$ and $q_{e,i}$ denote the active and reactive power generated by the generator at bus $i$, 
and $p_{l,i}$ and $q_{l,i}$ denote the active and reactive power consumed by the load at bus $i$. 
 
Let $
\mathbf{x} =
[\boldsymbol{\delta}^\top\;\boldsymbol{\omega}^\top]^\top
\in \mathbb{R}^{2n_g},
\mathbf{s} =
[\mathbf{v}^\top\;\boldsymbol{\theta}^\top]^\top
\in \mathbb{R}^{2n_b},
$ denote the dynamic and algebraic variables, respectively with $n_g$ generators and $n_b$ buses. Then, the overall system is described by the differential algebraic equation 

\begin{equation}
\dot{\mathbf{x}} = \mathbf{f}(\mathbf{x}, \mathbf{s}), 
\qquad
\mathbf{0} = \mathbf{g}(\mathbf{x}, \mathbf{s}),
\label{eq:dae}
\end{equation}

where $\mathbf{f} : \mathbb{R}^{2n_g} \times \mathbb{R}^{2n_b} \rightarrow \mathbb{R}^{2n_g}$ represents the right hand side terms  in \eqref{eq:swing1}-\eqref{eq:swing2} and $\mathbf{g} : \mathbb{R}^{2n_g} \times \mathbb{R}^{2n_b} \rightarrow \mathbb{R}^{2n_b}$ represents the algebraic constraints in \eqref{eq:active_balance}–\eqref{eq:reactive_balance}.



For the system measurements, we generate them through the relationship $
\mathbf{z}(t) = \mathbf{h}(\mathbf{x}(t), \mathbf{s}(t)) + \mathbf{a}(t),
\label{eq:measurement}
$ where $\mathbf{a}(t)$ denotes the attack signal and $
\mathbf{h} : \mathbb{R}^{2n_g} \times \mathbb{R}^{2n_b}
\rightarrow
\mathbb{R}^{m} $
maps the system states to the active and reactive power injections from the generators.

\vspace{-1em}
\section{Dynamic Attack Modeling}\label{attack_design}
Consider the power system partitioned into $d$ measurement areas, with the measurement vector $
\mathbf{z}(t) = [\mathbf{z}_1^\top(t),\dots,\mathbf{z}_d^\top(t)]^\top
$, where $\mathbf{z}_k(t) \in \mathbb{R}^{m_k}$ corresponds to the measurements in area $k$. 
Let $\mathcal{T} = \{t_0, t_1, \dots, t_N\}$ denote the discrete time instants over the attack planning horizon. 
The attacker is assumed to have access to the complete measurement time series 
$\{\mathbf{z}_k(t_i)\}_{t_i \in \mathcal{T}}$ for each area $k$.

At any time $t$, the attack injected in area $k$ is parameterized as $
\mathbf{a}_k(t) = c(t)\,\boldsymbol{\mu}_k,
\label{eq:attack_model}
$ where $c(t)$ is a known time-dependent attack policy, such as a sinusoidal function, and $\boldsymbol{\mu}_k \in \mathbb{R}^{m_k}$ defines the spatial pattern of the attack within area $k$. We assume here the attacker possess sufficient system knowledge to construct a
local linear approximation of the measurement model around the operating condition, inferred from historical measurements or leaked state
estimates. Under this assumption, the measurement residual is linearized around a nominal operating
point $(\mathbf{x}_0,\mathbf{s}_0)$ as:

\begin{equation}
h(\mathbf{x},\mathbf{s}) \approx h(\mathbf{x}_0,\mathbf{s}_0) + H\,\Delta \mathbf{s},
\end{equation}

where $\Delta \mathbf{s}= [\mathbf{x}^\top\;\mathbf{s}^\top]^\top
- [\mathbf{x}_0^\top\;\mathbf{s}_0^\top]^\top$ and $H$ is the measurement Jacobian. Given a WLS state estimator, the estimated state deviation is given by:

\begin{equation}\label{residual}
\Delta \hat{\mathbf{s}} = \arg \min_{\Delta \mathbf{s}} 
\|\mathbf{z} - H \Delta \mathbf{s}\|_W^2 
= (H^\top W H)^{-1} H^\top W \mathbf{z}.
\end{equation}

Substituting this solution into the attack induced part of the residual $\mathbf{r}(t) = \mathbf{z}(t) - H \Delta \hat{\mathbf{s}}(t)$ with $ \mathbf{z}(t) \approx  H \Delta s + \mathbf{a}(t)$ yields the approximation:

\begin{equation}
\mathbf{r}(t) \approx R\,\mathbf{a}(t), \qquad 
R = I - H(H^\top W H)^{-1} H^\top W,
\end{equation}

where $R \in \mathbb{R}^{m \times m}$. 
Consequently, the residual in
area $k$ is approximated as
$
\mathbf{r}_k(t) \approx c(t)\,R_k \boldsymbol{\mu}_k,
$ where $R_k \in \mathbb{R}^{m\times m_k}$. This residual is critical for measuring the stealthiness of the attacks. A bad data detector is bypassed when $\norm{\mathbf{r}_k(t)} \leq \epsilon_k$, where $\epsilon_k$ is the detection threshold. To measure the effectiveness of the attacks, we measure the deviation of the state estimate under attack in \eqref{residual} from the state estimate under nominal conditions. Accordingly, the attack effectiveness is given as:

\begin{align}
    \mathcal{E}(t) \approx \| Q \mathbf{a}(t)\|, \qquad  Q = (H^\top W H)^{-1} H^\top W,
\end{align}

where $Q \in \mathbb{R}^{n \times m}$, $n = 2n_g+2n_b$. And for each area $k$, it is approximated as $\mathcal{E}_k(t) \approx \| c(t)\,Q_k \boldsymbol{\mu}_k \|
$, where $Q_k \in \mathbb{R}^{n \times m_k}$.
Thus, the attack design problem is formulated as:
\vspace{-2em}
\begin{equation}
\begin{alignedat}{2}
\max_{\{\boldsymbol{\mu}_k\}_{k=1}^d} \quad
& \alpha
  \sum_{k=1}^{d} \|Q_k \boldsymbol{\mu}_k\|_2^2 \\
& {} - \sum_{k,l=1}^d w_{kl}
  \| Q_k\boldsymbol{\mu}_k - Q_l\boldsymbol{\mu}_l \|_2^2 \\
\text{s.t.} \quad
& \| c(t)R_k \boldsymbol{\mu}_k \|_2 \le \epsilon_k, \\
& \|\boldsymbol{\mu}_k\|_1 \le \rho_k.
\end{alignedat}
\end{equation}

where $\alpha = \sum_{t \in \mathcal{T}} c(t)^2$, $w_{kl}$ quantifies inter-area electrical coupling via impedances, and $\rho_k$ controls the sparsity level of the attack. The formulation is nonconvex due to the maximization of convex quadratic terms (since $Q_k^\top Q_k \succeq 0$). To obtain a tractable solution, we compute a locally optimal solution using the projected gradient ascent (PGA) method presented in Algorithm~\ref{projection_attack_design}.
\begin{algorithm} 
\caption{PGA for Attack Design}
\begin{algorithmic}[1]
\State Initialize $\{\boldsymbol{\mu}_k^{(0)}\}_{k=1}^d$
\For{$t = 0$ to $T-1$}
    \For{$k = 1$ to $d$}
        \State Compute gradient:

        \[
\begin{aligned}
\nabla_{\boldsymbol{\mu}_k} f
&= 2\alpha Q_k^\top Q_k \boldsymbol{\mu}_k^{(t)} \\
&\quad - 2 \sum_{l=1}^d w_{kl}
Q_k^\top \big( Q_k \boldsymbol{\mu}_k^{(t)} - Q_l \boldsymbol{\mu}_l^{(t)} \big)
\end{aligned}
\]

        \State Gradient step:
        $
        \tilde{\boldsymbol{\mu}}_k \leftarrow \boldsymbol{\mu}_k^{(t)} + \eta \nabla_{\boldsymbol{\mu}_k} f
        $
        
        \State Projection:

\[
\begin{aligned}
\boldsymbol{\mu}_k^{(t+1)} \leftarrow \arg\min_{\boldsymbol{\mu}_k} \quad
& \|\boldsymbol{\mu}_k - \tilde{\boldsymbol{\mu}}_k\|_2^2 \\
\text{s.t.} \quad
& \|c(t)R_k \boldsymbol{\mu}_k\|_2 \le \epsilon_k, \\
& \|\boldsymbol{\mu_k}\|_1 \le \rho_k.
\end{aligned}
\]

\EndFor
\EndFor
\State \Return $\{\boldsymbol{\mu}_k\}_{k=1}^d$
\end{algorithmic}\label{projection_attack_design}
\end{algorithm}

\section{Kernel-Embedded Detection}
The KEFSD framework detects coordinated attacks by learning the intrinsic structure of residual trajectories. Residuals are first modeled in an RKHS to capture nonlinear dynamics, and a nominal subspace is constructed using RKHS-constrained functional PCA. Online, residual trajectories are projected onto this subspace, and deviations are quantified for real-time anomaly detection.

\subsection{RKHS Modeling of Residual Trajectories}
Let $\mathbf{r}(t)\in\mathbb{R}^{m}$ denote the clean residual vector generated by the state estimator, where $m=\sum_{k=1}^{d} m_k$ is the total
number of residual channels across all $d$ areas. For each channel
$i\in\{1,\dots,m\}$, the residual trajectory is denoted by $r_i(t)$ and is
observed over the interval $t\in[0,T_s]$ where $T_s$  is chosen sufficiently large to ensure that the power system has reached steady-state operation under nominal conditions. The residuals are sampled at discrete time instants $\{t_j\}_{j=1}^{n}\subset[0,T_s]$.

To enable a kernel-based functional representation, each residual channel is approximated by a latent function
$g_i^{\text{nom}}(\cdot)\in\mathcal{H}_\kappa$. Specifically, we adopt the standard
observation model:

\begin{equation}\label{eq:rkhs_obs_model}
r_i(t_j)=g_i^{\text{nom}}(t_j)+w_i(t_j),\qquad j=1,\dots,n,
\end{equation}

where $w_i(t_j)$ captures measurement noise. Note that
\eqref{eq:rkhs_obs_model} does not require $r_i(\cdot)$ to
belong to $\mathcal{H}_\kappa$; rather, it assumes that the residual trajectory
can be well-approximated by an RKHS function.

By the representer theorem~\cite{scholkopf2001generalized}, the
RKHS function $g_i^{\text{nom}}(\cdot)$ admits the finite kernel expansion:

\begin{equation}\label{eq:representer}
g_i^{\text{nom}}(t)=\sum_{j=1}^{n}\alpha_{ij}\,\kappa(t,t_j),
\end{equation}


Let $K\in\mathbb{R}^{n\times n}$ denote the Gram matrix associated with the
kernel $\kappa$, with entries
$
K_{ij}=\kappa(t_i,t_j), i,j\in\{1,\dots,n\}
$. Evaluating \eqref{eq:representer} at the sampling instants yields

$\mathbf{g}_i^\text{nom} =
\begin{bmatrix}
g_i^\text{nom}(t_1) & \cdots & g_i^\text{nom}(t_n)
\end{bmatrix}^\top
=K\boldsymbol{\alpha}_i
$, where $\boldsymbol{\alpha}_i=[\alpha_{i1}\;\cdots\;\alpha_{in}]^\top$. 

Thus, the RKHS approximation is given as:

\begin{equation}\label{eq:rkhs_regression}
\hat g_i^{\text{nom}}
=
\arg\min_{g_i^{\text{nom}}\in\mathcal{H}_\kappa}
\sum_{j=1}^{n}\big(r_i(t_j)-g_i^{\text{nom}}(t_j)\big)^2
+\lambda \|g_i^{\text{nom}}\|_{\mathcal{H}_\kappa}^2,
\end{equation}

where $\lambda>0$ controls the trade-off between fidelity to the sampled
residuals and smoothness of the RKHS approximation. Substituting  \eqref{eq:representer} into \eqref{eq:rkhs_regression} yields:

\begin{equation}\label{eq:alpha_opt}
\hat{\boldsymbol{\alpha}}_i
=
\arg\min_{\boldsymbol{\alpha}\in\mathbb{R}^{n}}
\|\mathbf{r}_i-K\boldsymbol{\alpha}\|_2^2
+\lambda\,\boldsymbol{\alpha}^\top K\boldsymbol{\alpha},
\end{equation}

with closed-form solution, $
\hat{\boldsymbol{\alpha}}_i=(K+\lambda I_n )^{-1}\mathbf{r}_i$, where $
\mathbf{r}_i =
\begin{bmatrix}
r_i(t_1) & r_i(t_2) & \cdots & r_i(t_n)
\end{bmatrix}^\top \in\mathbb{R}^{n}
$. Consequently, the estimated latent residual function is:

\begin{equation}\label{eq:g_hat}
\hat g_i^{\text{nom}}(t)=\sum_{j=1}^{n}\hat\alpha_{ij}\,\kappa(t,t_j).
\end{equation}

\subsection{RKHS-Constrained Functional PCA (fPCA)} \label{functional_pca}
Let $\mathbf{y}_i \in \mathbb{R}^n$, $i = 1,\dots,m$, denote noisy discrete observations of $g_i^{\text{nom}} : [0,T_s] \rightarrow \mathbb{R}$ evaluated at sampling points $t_j \in [0,T_s]$, $j = 1,\dots,n$. The goal here is to identify an $r$-dimensional functional subspace:
$
\mathcal{F}^* = \mathrm{span}\{f_1,\dots,f_r\} \subset \mathcal{H}_\kappa
$ that captures the dominant variation of the latent functions. 
 Now to capture the dominant variability in the observed residuals, we examine the covariance matrix:

\begin{equation} \Sigma_m = \frac{1}{m}  Y Y^\top, \end{equation}

where $Y = [\, \mathbf{y}_1 \mid \mathbf{y}_2 \mid \dots \mid \mathbf{y}_m \,]$.
In classical PCA, principal directions are obtained from the eigendecomposition of $\Sigma_m$. However, such directions are unconstrained and may correspond to non-smooth functions. To enforce smoothness of the extracted functional components, we incorporate an RKHS norm constraint.

For any vector $z \in \mathbb{R}^n$, we define the associated RKHS function:

\begin{equation}
f(t) = \sum_{j=1}^{n} (K^{-1} z)_j \kappa(t_j,t).
\end{equation}

From RKHS theory \cite{scholkopf2001generalized}, its squared norm satisfies $\|f\|_{\mathcal{H}_\kappa}^2 = z^\top K^{-1} z$. Therefore, controlling the quadratic form $z^\top K^{-1} z$ is equivalent to controlling the RKHS norm of the corresponding functional direction.

Motivated by this observation, we consider the constrained optimization problem:




\begin{equation}
\begin{aligned}
\max_{Z \in \mathbb{R}^{n \times r}} \quad 
& \mathrm{Tr}\!\left(Z^\top (\Sigma_m - \gamma K^{-1}) Z\right) \\
\text{s.t.} \quad 
& Z^\top Z = I_r.
\end{aligned}\label{kpca_opt}
\end{equation}

where the optimal matrix $Z$ consists of the eigenvectors corresponding to the largest $r$ eigenvalues. In practice, $\gamma$ is chosen from a grid $\gamma \in [\gamma_{\min}, \gamma_{\max}]$ and for each candidate value, $\|f\|_{\mathcal{H}_\kappa} = \left\langle Z_\gamma, K^{-1} Z_\gamma \right\rangle $ is evaluated, where $Z_\gamma$ denotes the optimizer of \eqref{kpca_opt}.  The value of $\gamma$ yielding the minimum admissible RKHS norm is then selected. On close inspection it can be seen the solution to the problem is actually the eigendecomposition of $\Sigma_m - \gamma K^{-1}$. Finally, the functional principal components are reconstructed as:

\begin{equation}
f_i(t) =
\sum_{j=1}^{n}
\left(K^{-1} Z[i]\right)_j
\kappa(t_j,t),
\quad i=1,\dots,r.
\end{equation}

The process described above is depicted in Algorithm~\ref{RKHS_constrained}.

 \vspace{-1em}
\begin{algorithm}
\caption{RKHS-Constrained $f$PCA}
\begin{algorithmic}[1]
\Procedure{$f$PCA}{$Y, \kappa, r$}
    \State $K_{ij} \leftarrow \kappa(t_i, t_j), \; t_i \in [0,T_s]$ \Comment{Gram matrix}
    \State $\Sigma_m \leftarrow \mathrm{cov}(Y)$ \Comment{Sample covariance matrix}
    \State $ \boldsymbol{\gamma} \leftarrow \mathrm{linspace}(\gamma_{\min}, \gamma_{\max}, n)$
    \For{$i = 1$ to $n$}
        \State $X \leftarrow \mathrm{eig}(\Sigma_m - \gamma_i K^{-1})$
        \State $Z[i] \leftarrow X_{1:r}$
        \Comment{assign top $r$ eigenvalues}
        \State $T[i] \leftarrow \langle Z[i], K^{-1} Z[i] \rangle$
    \EndFor
    \State $i^* \leftarrow \arg\min_i T[i]$ \Comment{$i \in [1,n]$}
    \State \Return $Z[i^*]$
\EndProcedure
\end{algorithmic}\label{RKHS_constrained}
\end{algorithm}

\subsection{Online Projection-Based Attack Detection}
To monitor residuals in real time, we adopt a sliding-window strategy. 
Specifically, at each sampling instant $t_s$ and for each channel $i \in \{1,\dots,m\}$, 
the most recent window: $\mathcal{W}_s = \{t_{s-W+1}, \dots, t_s\}$ is used to estimate a latent residual trajectory $\hat g_i \in \mathcal{H}_\kappa$ via the regularized RKHS problem:

\begin{equation}
\hat{g}_i
=
\arg\min_{g_i \in \mathcal{H}_\kappa}
\sum_{j \in \mathcal{W}_s}
\big(r_i(t_j) - g_i(t_j)\big)^2
+ \lambda \|g_i\|_{\mathcal{H}_\kappa}^2,
\end{equation}

where $\lambda > 0$ controls the smoothness–fidelity tradeoff.

Although $\hat g_i \in \mathcal{H}_\kappa$, it may deviate from the nominal functional subspace 
$\mathcal{F}^* = \mathrm{span}\{f_1,\dots,f_r\}$ learned from attack-free data. 
Its orthogonal projection onto $\mathcal{F}^*$ is therefore:

\begin{equation}
\Pi_{\mathcal{F}^*} \hat g_i
=
\sum_{\ell=1}^{r}
\langle \hat g_i, f_\ell \rangle_{\mathcal{H}_\kappa} \, f_\ell,
\end{equation}

and the resulting projection residual energy is:

\begin{equation}
J_i(t_s)
=
\left\|
(I-\Pi_{\mathcal{F}^*}) \hat g_i
\right\|_{\mathcal{H}_\kappa}^2.
\end{equation}

To obtain an area-level detection statistic, the channel-wise energies are aggregated so that, for area $k$ containing $m_k$ channels, we have:

\begin{equation}
J^{(k)}(t_s)
=
\frac{1}{m_k}
\sum_{i=1}^{m_k}
\left\|
(I-\Pi_{\mathcal{F}^*}) \hat g_i
\right\|_{\mathcal{H}_\kappa}^2.
\end{equation}

Consequently, an alarm is triggered if $J^{(k)}(t_s) > \tau^{(k)}$, where $\tau^{(k)}$ is the detection threshold for each area.

\section{Simulation Experiments}
\begin{figure}[htpb]
    \centering
    \includegraphics[scale=0.6]{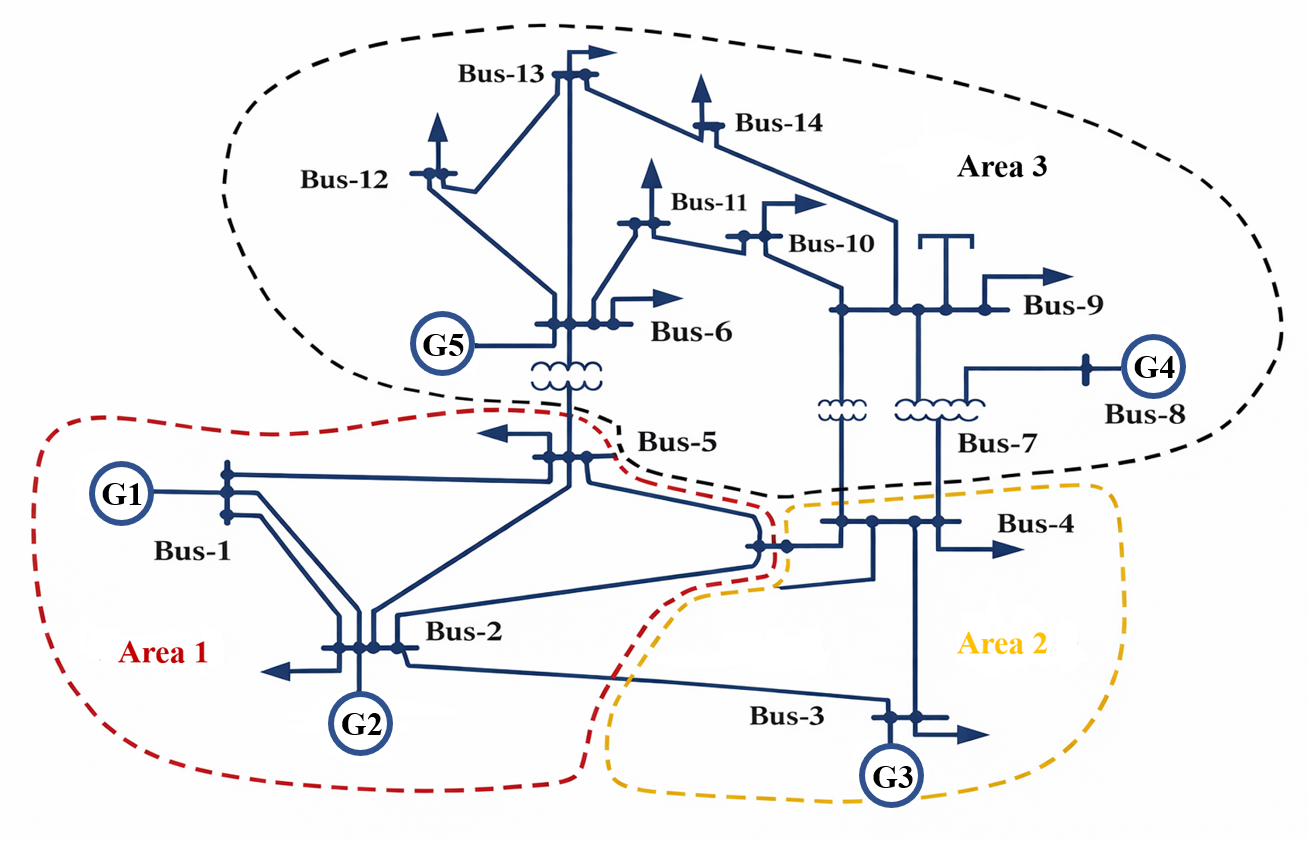}
    \caption{Modified IEEE 14-bus power network divided into three areas (Area 1–3), highlighting generators (G1 - G5) and interconnecting transmission lines.}
    \label{IEEE_14Bus_new}
\end{figure}

The simulation studies are carried out on the IEEE 14-bus power system, which has been partitioned into three distinct areas ($k=3$) as illustrated in Fig.~\ref{IEEE_14Bus_new}. In this configuration, Generators 1 and 2 are located in Area 1, Generator 3 in Area 2, and Generators 4 and 5 in Area 3. The dynamic behavior of each generator is captured using the classical swing equation (see Section \ref{dynamic_power_model}), with inertia constants, $m = [2.1, 2.2, 2.3, 2.4, 2.5]$~s and damping coefficients, $d = [0.72, 0.71, 0.73, 0.65, 0.70]$~p.u., chosen within typical ranges for synchronous generators~\cite{kundur2007power}. The matrices, $G$ and $B$ are obtained from the system topology and line parameters. The operating point $(\mathbf{x}_0, \mathbf{s}_0)$ is taken as the steady-state IEEE 14-bus solution, where $\mathbf{x}_0$ consists of the generator rotor angles and frequency ($120\pi$ rad/s), and $\mathbf{s}_0$ consists of the bus voltage magnitudes and phase angles from the IEEE bus data. To account for interactions between the different areas, the inter-area weighted adjacency matrix, $[w_{kl}]$ is derived by considering the impedance-based coupling between connected areas.

\begin{figure}[htbp]
    \centering
    \includegraphics[scale=1]{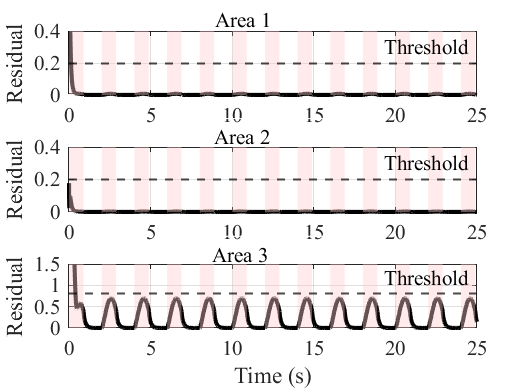}
    \caption{Residual norm trajectories with shaded attack intervals ($45 \%$ duty cycle) across areas. The residual stays below the threshold, indicating successful stealthy attacks.}
    \label{residual_norm_signals}
\end{figure}

The system outputs consist of $m=10$ measurements corresponding to the active and reactive power of the generators. For the nominal model, the detection threshold is set to $\epsilon = [0.2, 0.2, 0.8]$. Attack scenarios, as formulated in Section~\ref{attack_design}, are generated using $\boldsymbol{\mu}_k$ and $c(t)$ chosen as a sinusoidal signal of frequency 1~Hz, zeroed at a 45\% duty cycle. The resulting residual norm trajectories are presented in Fig.~\ref{residual_norm_signals}. As seen, the residual remains below the predefined thresholds, demonstrating the stealthy nature of the injected attacks. To learn the nominal subspace $\mathcal{F}^*$, $n=200$ clean residual samples were collected at $t_s=0.01$s. The latent residual function $\hat{g}^{\text{nom}}$ was estimated using a Gaussian kernel with regularization parameter $\lambda=10^{-3}$. The principal components were then obtained via Algorithm~1, using $Y$ formed by stacking the estimated latent residuals, a Gaussian kernel $\kappa$, and subspace dimension $r=200$; the parameter $\gamma$ was swept over $[10^{-6}, 10]$. During real-time operation, residual trajectories were estimated using a sliding window of size $W=20$. Detection thresholds were tuned to match the false alarm rate of the residual norm-based detector, resulting in $\tau = [0.45, 0.56, 0.78]$.

\begin{figure}[t]
    \centering
    \begin{subfigure}[b]{0.45\textwidth} 
        \centering
        \includegraphics[scale=1]{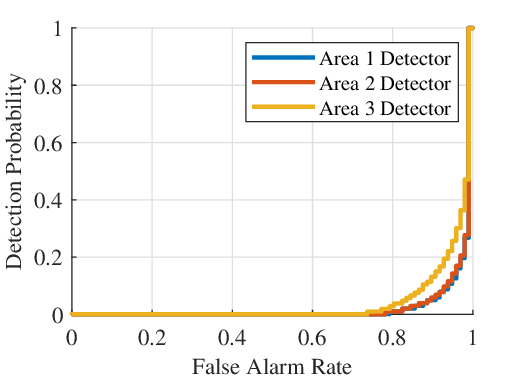}
        \caption{}
        \label{fig:detector1}
    \end{subfigure}
    \hspace{0.5\textwidth} 
    \begin{subfigure}[b]{0.45\textwidth} 
        \centering
        \includegraphics[scale=1]{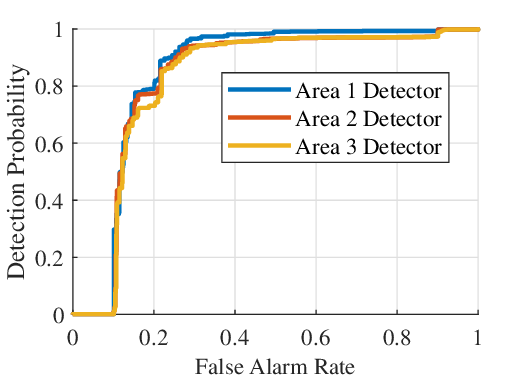}
        \caption{}
        \label{fig:detector2}
    \end{subfigure}

    \caption{Detection performance across bus network areas: (a) residual norm–based detector with delayed response, (b)proposed  KEFSD-based detector with faster and more reliable detection.}
    \label{fig:detectors_comparison}
\end{figure}

Table 1 shows the performance scores of the detector such as true positive rates (TPR), false positive rates (FPR), false negative rates (FNR), precision, and F1 score. As seen in Table \ref{Metrics}, the detector achieves high TPR (78–83\%) and precision (88–90\%) while FNR (17–22\%). The F1 scores above 86\% indicate a strong balance between detection sensitivity and reliability, demonstrating the effectiveness of KEFSD in identifying attacks across all areas. In Fig. \ref{fig:detectors_comparison}, it is clear that the proposed KEFSD detector significantly outperforms the residual norm–based method, achieving high detection probability at much lower false alarm rates. Table \ref{tab:five_col} further quantifies this improvement using the area under the detection probability–false alarm rate curves (AUC), which reflects the detector’s ability to distinguish between normal and attack conditions. The consistently higher AUC values of KEFSD across all areas confirm its superior overall detection performance.

\begin{table}[htbp]
\caption{Performance metrics of the KEFSD across areas.}
\centering
\begin{tabular}{|c|c|c|c|}
\hline
Metrics & Area 1 (\%)  & Area 2 (\%) & Area 3 (\%) \\ 
\hline
TPR &78.30 &81.10  &83.11  \\ 
\hline
FPR &13.33  &12.23  &10.23   \\ 
\hline
FNR &21.7  &18.90  &16.89    \\ 
\hline
Precision &88.34 &89.20 &89.56    \\ 
\hline
F1 score &87.45  &86.78 &90.23  \\ 
\hline
\hline
\end{tabular}

\label{Metrics}
\end{table}
\vspace{-0.5em}
 \begin{table}[htbp]
\caption{AUC comparison of detectors across areas.  }
\centering
\begin{tabular}{|c|c|c|c|}
\hline
Detectors& Area 1 (\%) & Area 2 (\%) & Area 3 (\%) \\ 
\hline
Residual norm  &74.5  &77.4  &75.8  \\ 
\hline
KEFSD  &84.80  &83.0  &82.4    \\ 
\hline
\hline
\end{tabular}
\label{tab:five_col}
\end{table}

\section{Conclusion}
This paper presented a joint framework for modeling and detecting coordinated dynamic stealthy attacks in power systems. A time-aggregated multi-area attack model was developed to capture temporal evolution and inter-area coordination, along with a tractable attack design formulation. For detection, a KEFSD method was proposed, leveraging RKHS-constrained functional PCA to identify structured temporal deviations that missed by residual-based methods. Simulation results on the IEEE 14-bus system show that KEFSD achieves improved detection performance and robustness under coordinated stealthy attacks. These results highlight the importance of incorporating temporal and functional structure in power system security. Future work includes extending the framework to large-scale systems, and developing adaptive and distributed implementations.

\section*{Code Availability}
The source code used to generate the simulation results presented in this paper
is made publicly available at:
https://github.com/resilient-autonomous-systems-lab.
\section*{Acknowledgment}
This work was supported by the U.S. Department of Energy
under Award Number DE-CR0000028. The authors
acknowledge that the views expressed do not necessarily
reflect those of the United States Government.

\bibliographystyle{ieeetr}

\bibliography{myreferences}             
\end{document}